\documentclass[twoside,openright,11pt,a4paper,epsf]{article}
\usepackage{latexsym}
\usepackage[T1]{fontenc}
\usepackage{amsmath}
\usepackage{amsthm}
\usepackage{amsfonts} 
\usepackage{amssymb}
\usepackage{vmargin}
\usepackage{epsfig}
\usepackage[all]{xy}

\newcommand{\color}[6]{}

\newcommand{\N}{\mathbb{N}}

\newcommand{\C}{\mathbb{C}}

 \newcommand{\nbd}{neighbourhood }
\newcommand{\fonction}[5]
{$$ 
\begin{array}{rcccl}
 #1 & : & #2 & \longrightarrow &#3 \\
    &   & #4 & \longmapsto &#5 
\end{array}
$$}

\newcommand{\spc}{\mathcal{SPC}}

\newcommand{\re}{\text{Re}\,}

\newcommand{\grad}{\vec{\nabla}}

\newcommand{\priv}{\backslash}
\newcommand{\lra}{\longrightarrow}

\newcommand{\om}{\Omega}
\newcommand{\eps}{\varepsilon}
\renewcommand{\phi}{\varphi}

\renewcommand{\L}{\mathcal{L}}
\newcommand{\psh}{$p.s.h$ }

\newcommand{\wdt}[1]{\widetilde{#1}}
\newcommand{\cqfd}{\hfill $\square$ \vspace{0.1cm}\\ }
\newcommand{\sbull}{{\tiny $\bullet$}}

\newtheorem*{thm*}{Theorem}
\newtheorem{definition}{Definition}[section]
\newtheorem{lemma}[definition]{Lemma}
\newtheorem{thm}{Theorem}

\newtheorem{cor}[definition]{Corollary}
\newtheorem{prop}[definition]{Proposition}

\title{\vspace{-1.5cm} A Wong-Rosay type theorem for proper holomorphic 
self-maps.}
\author{Emmanuel Opshtein}\date{}
\begin{document}
\maketitle
\begin{abstract}
In this short paper, we show that the only proper holomorphic self-maps of bounded 
domains in $\C^k$ whose dynamics escape to a strictly pseudoconvex point of 
the boundary are automorphisms of the euclidean ball. This is a Wong-Rosay type theorem for  
a sequence of maps whose degrees are {\it a priori} unbounded.
\end{abstract}

\section*{Introduction.}
In 1977, Wong proved that the only strictly pseudoconvex domain with non-compact
automorphism group is the ball \cite{wong}. This result was generalized
by Rosay \cite{rosay} (see also \cite{pinchuk}).
\begin{thm*}[Wong-Rosay]
Let $\om$ be a bounded domain in $\C^k$ and $(f_n)$ a sequence of its automorphisms.
 Assume that the orbit of a point of $\om$ under $(f_n)$ accumulates a smooth 
strictly pseudoconvex point of $b\om$.
Then $\om$ is biholomorphic to the euclidean ball.
\end{thm*}
This theorem remains valid for a sequence of correspondences provided that their degrees 
remain bounded \cite{ourimi2}. In this paper, we  prove that the theorem above also holds true in presence of unbounded degree, when the sequence of automorphisms is replaced by the iterates of a proper 
holomorphic self-map.
\begin{thm}\label{properwong}
Let $\om$ be a bounded domain in $\C^k$ with a proper holomorphic self-map $f$. 
If there is a point $y$ of  $\om$ whose orbit under the iterates of $f$ accumulates a smooth strictly 
pseudoconvex point $a$ of $b\om$ (that is $f^{n_k}(y)\lra a$), 
then $\om$ is biholomorphic to the euclidean ball
and $f$ is an automorphism.
\end{thm}
In \cite{moi1}, the question of wether a proper holomorphic self-map of a smoothly bounded domain in $\C^k$ 
has to be an automorphism of the domain was considered. In $\C^2$ for instance, it was proved that non-injective proper self-maps of such domains has a non-compact dynamics (all the limit maps of the dynamics have value on the 
boundary of $\om$). Theorem \ref{properwong} goes one step further in this direction : the limit maps even take values in the weakly pseudoconvex part of the boundary.

The main ingedient for this result is a local
version of Wong-Rosay's theorem concerning sequences of CR-maps. It 
was first obtained by Webster in the wake of Chern-Moser's theory of strictly pseudoconvex
hypersurfaces \cite{webster}.

\begin{thm}[Webster]\label{wrlocal} Let $(\Sigma,a)$ and $S$ be two germs of strictly pseudoconvex 
hypersurfaces. Assume there is a sequence of CR-embeddings of $S$ into 
$\Sigma$ whose images converge to $a$. Then $S$ is spherical, i.e. locally
CR-diffeomorphic to the euclidean sphere.
\end{thm}
The idea behind the proof of theorem \ref{properwong} is to consider the CR-maps 
induced by $f^n$ on the boundary rather than the maps $f^n$ themselves. Using techniques developped 
in \cite{moi1}, we study the way these CR-maps degenerate and check that theorem \ref{wrlocal} 
applies : around $a$, the boundary of $\om$ is spherical. The local biholomorphism between our domain and the ball then extends 
to the whole of $\om$ thanks to the dynamical situation.  \vspace{.2cm}

The paper is organised as follows. We first collect some trivial dynamical facts about $f$ and the automorphisms of the ball which will allow us to propagate the local sphericity to the whole domain.
In section \ref{proof}, we prove theorem \ref{properwong} {\it  modulo} the central question of the local sphericity around $a$. In section \ref{locspher}, we finally turn back to this problem.

\section{Preliminary remarks.}
Surprisingly enough, the convergence hypothesis on $f^{n_k}(y_0)$ in theorem \ref{properwong} has very strong (though very classical) implications in the holomorphic context. The aim of this section is to clarify some of them, as well as pointing out the well-known properties of the dynamics of the automorphisms of the ball which will be usefull to us. Henceforth, $\om$, $f$ and $a$ are as in theorem \ref{properwong}. \vspace{.2cm}\\
\indent{} First of all, this hypothesis may seem weaker than it actually is. Indeed, the contracting property of holomorphic maps for the Kobayashi distance (which is a genuine distance on bounded domains) leads to the following classical fact :
\begin{lemma}\label{oneall}
Any sequence of holomorphic maps between bounded domains $\om$ and $\om'$, which takes a point $y$ in $\om$ to a sequence converging to a strictly pseudoconvex point of the boundary of $\om'$ converges locally uniformly to this point on $\om$. For instance, the sequence $f^{n_k}$ converges locally uniformly to $a$ on $\om$.
\end{lemma}
\begin{cor}\label{sext}
The map $f$ extends smoothly to a \nbd of $a$ in $b\om$ and $f(a)=a$. Moreover, $f$ is a local biholomorphism 
(resp. CR-automorphism) in a \nbd of $a$ (resp. in $b\om$).
\end{cor}
\noindent {\it Proof :} Call $z_k:=f^{n_k}(y)$ and $w_k:=f(z_k)$. Since $w_k=f(f^{n_k}(y))=f^{n_k}(f(y))$,
 both $z_k$ and $w_k$ tend to $a$ because of the previous lemma. Since $a$ is a strictly pseudoconvex point, an observation of Berteloot ensures that $f$ extends continuously to a \nbd of $a$ in $b\om$  \cite{berteloot4}, with $f(a)=\lim f(z_k)=\lim w_k=a$. Such an extension is automatically smooth because $a$ is a strictly pseudoconvex smooth point of $b\om$  \cite{bell}. Since we are close to a strictly pseudoconvex point of the boundary, branching is prohibited and $f$ must be one-to-one (see \cite{fornaess}). \cqfd

Let us now discuss the dynamical type of the fixed point $a$. Although it attracts part of the dynamics, it is not obvious at first glance that $a$ is not a repulsive fixed point. The orbit of $y_0$ could {\it in principle} jump close to $a$ from time to time, then get expelled away from $a$. The following lemma shows that such a behaviour does not occur in our holomorphic context.
\begin{lemma} \label{nr}
The point $a$ is a non-repulsive fixed point of $f$. 
\end{lemma}
\noindent {\it Proof :} Assume by contradiction that $f$ is repulsive at $a$. By definition, 
there is an open \nbd $U$ of $a$ on which the inverse $f^{-1}$ of $f$ is well defined, takes 
values in $U$, and is even contracting : $d(f^{-1}_{|U}(z),a)<d(z,a)$ for any $z\in U$.
By assumption, there is a point $y_0\in \om$ such that 
$f^{n_k}(y_0)\in U$ as soon as $k$ is large enough. Define then
$$
n_k':=\min\{n\;|\; f^i(y_0)\in U, \hspace{.5cm} \forall i\in [n,n_k] \},
$$
so that $f^{n_k'-1}(y_0)\notin U$.
Since $f^{-1}_{|U}$ is contracting, the point $f^{n_k'}(y_0)$ is closer 
to $a$ than $f^{n_k}(y_0)$, so it tends to $a$ (in particular $(n_k')$ is an 
extraction). Equivalently $f^{n_k'-1}(f(y_0))$ tends to $a$, so $f^{n_k'-1}$ 
converges locally uniformly to $a$ by lemma \ref{oneall}. This is 
in contradiction with $f^{n_k'-1}(y_0)\notin U$.\cqfd 
\paragraph{} Let us finally discuss the dynamics of the automorphisms of the ball. Since there are very few of them (they form a finite dimensional group), their dynamics is rather poor and any small piece of information on it may give rise to strong restrictions.
Recall the following well-known classification (see \cite{rudin}, section 2.4).
\begin{prop}\label{class}
Let $g$ be an automorphism of the unit ball in $\C^n$. Then the dynamics of $g$ is 
\begin{itemize}
\item[\sbull] either hyperbolic (North-South) : there exist exactly two fixed points $N,S\in bB$
of $g$ and $g^n$ converges locally uniformly to $S$ on $\overline{B}\priv\{N\}$.
\item[\sbull] or parabolic (South-South) : there exists a unique fixed point $S\in bB$ of $g$
and $g^n$ converges locally uniformly to $S$ on $\overline{B}\priv\{S\}$.
\item[\sbull] or recurrent (compact) : The $g$-orbits remain at fixed distance from $bB$. 
If $g$ has a fixed point on $bB$ then it has a whole complex pointwise fixed line through this point
(see also \cite{mccluer}).
\end{itemize}
\end{prop}
\noindent What will be of interest for us in this classification is contained in the following lemma,
whose proof is straightforward from the classification.
\begin{lemma}\label{dynrest}
Let $g$ be a ball automorphism which has a non-repulsive fixed point $p$ on $bB$, and no interior fixed point near $p$. Then the dynamics of $g$ is either hyperbolic or parabolic, with south pole $p$ (meaning that $S$ is $p$ in the previous classification). Moreover, given any \nbd $U$ of $p$, there is a point $z$ in $U$ whose orbit remains in $U$ and converges to $p$. 
\end{lemma}

\section{Proof of theorem \ref{properwong}.}\label{proof}
In this section, we prove theorem \ref{properwong} leaving aside the central question of the sphericity of 
$b\om$ around $a$, which will be dealt with in the next section.  Let us first fix the notation. Let
$(\om,f,a)$ be a triple as in theorem \ref{properwong}. By a global change of coordinates in $\C^k$, we can 
take $a$ to the origin, the tangent plane of $b\om$ at $a$ to $\{\re z_1=0\}$, and make $\om$ strictly convex 
locally near $a$. For $\alpha$ small enough, define $U_\alpha$ and
$\om_\alpha$ as being the connected components of $a$ in $b\om\cap
\{\re z_1<\alpha\}$ and $\om\cap \{\re z_1<\alpha\}$. \vspace{,2cm}

The first step of the proof, postponed to the following section, consists in showing that $b\om$ is spherical around $a$. 
\begin{lemma}\label{spher}
A neighbourhood of $a$ in $b\om$ is spherical.
\end{lemma}
\noindent This means that there exists a CR-diffeomorphism $\Phi:U_\eps\lra V\subset bB$. A classical extension theorem even shows that $\Phi$ extends to a biholomorphism $\Phi:\om_\eps\lra D$ where $D$ 
is an open set of $B$ whose boundary contains $V$ (see \cite{boggess}). This biholomorphism allows 
to transport $f$ to a local automorphism of $B$, defined by
\fonction{g}{\Phi(\om_\eps\cap f^{-1}(\om_\eps))}{\Phi(\om_\eps)}{x}{
\Phi\circ f\circ \Phi^{-1}(x).}  
The key point of the whole proof is the following extension
phenomena discovered by Alexander \cite{alexander} (see also \cite{pinchuk,hetu} for 
the form of the result we use here). The local biholomorphism $g$ uniquely extends
to a global automorphism of the ball, again denoted by $g$. \vspace{,2cm}

The second step consists in using the dynamics of $f$ and the injectivity of $g$ (which we got for free thanks to Alexander's theorem) to propagate the local sphericity, and produce a biholomorphism between $\om$ and $B$. 
Let us first discuss the possible dynamics of $g$. By lemma \ref{nr}, $a$ is not a repulsive fixed point of $f$
so  $\Phi(a)$ is neither one for $g$. Moreover, since  $f$ has no fixed point inside $\om$ (because of lemma  \ref{oneall}), $g$ has also no fixed point in $V$. By lemma \ref{dynrest}, $g$ is either hyperbolic with attractive fixed point $\Phi(a)$ or parabolic with only fixed point $\Phi(a)$. From now on, we will denote $S:=\Phi(a)$. The same  lemma  also guarantees that there are points in $D$ whose (positive) orbits under $g$ remain in $D$ and tend to $S$.  Since their whole orbits remain in $D$, the conjugacy thus allows to get the following informations on $f$ in return.
\begin{lemma}\label{cv}
The whole sequence of iterates  $(f^n)$ (rather than only a subsequence) converges to $a$ on $\om$. Moreover, the set
$$
\om_\eps':=\{z\in \om_\eps\; , \;  f^n(z)\in\om_\eps\;\; \forall n\in\N\}
$$
is a non-empty open invariant set of $f$.
\end{lemma}
\noindent {\it Proof :} From the discussion above, we conclude that there is a point $y$ in $\om_\eps$ such that $f^n(y)$ remains in $\om_\eps$ (thus $\om_\eps'$ is not empty). Its orbit also converges to $a$. By lemma \ref{oneall}, $f^n$ must therefore converge to $a$ locally uniformly on $\om$. The set $\om_\eps'$ is obviously invariant by $f$. Finally, it is open because the Kobayashi metric decreases under $f$. Indeed, if $z$ is in 
$\om_{\eps}'$, so is a Kobayashi $\delta$-neighbourhood of this point (take $\delta:= d_K\big(\om\cap \{\re z_1=\eps\}\,,\,\text{Orbit}(z)\big)$). \cqfd
\begin{cor}\label{extension}
The map $\Phi$ extends to a holomorphic map from $\om$ to $B$.
\end{cor}
\noindent {\it Proof :} Let $O_i$ denote $f^{-i}(\om_\eps')$. Because of the invariance of $\om_\eps'$ 
 by $f$ and since $f^n$ converges to $a$ on $\om$, we conclude that $(O_i)_i$ is a growing sequence 
 of open sets which exhausts $\om$. Define therefore
\fonction{\Phi}{\om=\cup O_i}{B}{z\in O_i}{g^{-i}\circ \Phi_{|\om_{\eps}'}\circ f^{i}(z).}
This map is obviously holomorphic (because $\om_{\eps}'$ is open), and coincides with $\Phi$ on 
$\om_{\eps}'$. It is therefore an extension of $\Phi$ itself. \cqfd 
The remaining point to prove is that $\Phi$ is a biholomorphism. Let us first prove that it is proper.
\begin{lemma} The map $\Phi$ is a proper map from $\om$ to $B$.
\end{lemma}
\noindent {\it Proof :} Recall that the dynamics of $g$ is either hyperbolic or parabolic. Moreover,  
$\Phi^{-1}\circ g(w) = f\circ \Phi^{-1}(w)$ for any $w\in D$ such that $g(w)$ belongs to $D$ (recall that 
$\Phi:\om_\eps\lra D$ is a biholomorphism).  
A basic consequence of these two facts is that $\Phi(O_n\priv O_{n-1})$ goes to $bB$ with $n$. 
Indeed, the $f$-orbit of a preimage by $\Phi$ of a point $w$ in this set reaches $O_0=\om_\eps'$ 
only at time $n$, so the $g$-orbit of $w$ cannot remain in $D$ before the same time (if $g^k(w)\in D$ for $k\geq N$, then $f^k(\Phi^{-1}(w))=\Phi^{-1}(g^k(w))$ is in $\om_\eps$ for $k\geq N$ also). 
If $n$ is large, $\Phi(z)$ has to be very close to some pole of the dynamics which is either $S$ if $g$ is parabolic
or another point of $bB$ if $g$ is hyperbolic. Anyway $\Phi(z)$ is close to the boundary of $B$. 

For an arbitrary sequence  $(z_i)_{i\in\N}\in \om$ converging to $b\om$, we must 
show that $\Phi(z_i)$ tends to the boundary of $B$. For this, fix
a positive real number $\delta$ and  an integer $n_0$ such that 
$d(\Phi(O_n\priv O_{n-1}),bB)\leq \delta$ for all $n\geq n_0$, meaning that $\Phi(\om\priv O_{n_0})$ is $\delta$-close from $bB$. Split then $(z_i)$
into two subsequences, one containing all the elements which belong to $O_{n_0}$, the other one  
those which escape from $O_{n_0}$ :
$$
\begin{array}{rcl}
(z_i^1) &:= &\{z_{n_i}\in \{z_n\}\;|\; z_{n_i}\notin O_{n_0}\},\\
(z_i^2) &:= &\{z_{n_i}\in \{z_n\}\;|\; z_{n_i}\in O_{n_0}\}.
\end{array}
$$   
By construction, $d(\Phi(z_i^1),bB)\leq \delta$. Since $f^{n_0}(z_i^2)\subset O_0$ 
and since  $f^{n_0}$, $\Phi_{|O_0}$  and $g$ are proper maps,  
$\Phi(z_i^2)=g^{-n_0}\circ\Phi_{|O_0}\circ f^{n_0}(z_i^2)$ is also $\delta$-close 
to $bB$ for $i$ large enough.\cqfd
Finally, we need to show that $\Phi$ is a biholomorphism. It is not 
yet clear since there exist holomorphic coverings of the ball. Anyway 
we know that any proper map to a bounded domain has a finite degree 
(see \cite{rudin}, chap. 15). In particular, there 
is an integer $d$ which bounds the numbers of preimages of $\Phi$ :
$$
\#\Phi^{-1}(z)\leq d, \hspace{2cm} \forall z\in B.
$$
Notice now that the degree of $\Phi$ bounds this of $f^n$ for all $n$
because $\Phi=g^{-n}\circ\Phi\circ f^n$. The degree of $f^n$ is thus bounded
on one hand and equal to $($deg$f)^n$ on the other. So $f$ is an automorphism 
of $\om$. The injectivity of $\Phi$ is now immediate since 
$\Phi_{|O_i}=g^{-i}\circ \Phi_{|O_0}\circ f^i$ is a composition of injective 
maps for all fixed $i$.\cqfd

\section{Local sphericity near the attractive point.}\label{locspher}
In this last section, we prove lemma \ref{spher}, namely that a neighbourhood of $a$ in $b\om$ is spherical.
 We recall that all the results proved in the previous section used this fact, so we have to go back to the general situation of theorem \ref{properwong}. Nevertheless, remind that we can speak of the action of $f$ on $b\om$, at least close to $a$, thanks to lemma \ref{sext}. The idea behind this technical part
of the proof is based on previous results concerning behaviours of sequences of CR-maps 
(see \cite{moi1,moi2}). Unformally speaking, they explain that 
non-equicontinuous sequences of CR-maps on  strictly pseudoconvex
hypersurfaces dilate a certain (anisotropic) distance. The proof of the sphericity 
then goes as follows. Either $f^{n_k}$ converges to $a$ on $\spc(b\om)$ and theorem \ref{wrlocal} gives the sphericity.
Or $f^{n_k}$ is not equicontinuous on $\spc(b\om)$ and it is dilating.
Then the inverse branches of $f^{n_k}$ are contracting CR-diffeomorphisms and 
theorem \ref{wrlocal} gives the sphericity. 
Let us first fix the easy situation where $f^{n_k}$ converges to $a$ on $\spc(b\om)$.
\begin{prop}\label{easyspher}
Assume $f^{n_k}$ converges locally uniformly to $a$ on a neighbourhood of $a$ in $b\om$. 
Then $b\om$ is spherical near $a$.
\end{prop}
\noindent{\it Proof :} Theorem \ref{wrlocal} explains that it is enough to 
find a contracting sequence of CR-automorphisms on a neighbourhood 
of $a$. We are assuming here that $(f^{n_k})$ is a sequence of contracting 
CR-maps on a piece of $\spc(b\om)$. Also, corollary \ref{sext} shows that $f$ is a local 
diffeomorphism at $a$. We thus only need to prove that there 
is a fixed neighbourhood of $a$ on which all $f^{n_k}$ are injective. To see this, 
first assume that $f^n$, and not only $f^{n_k}$, converges to $a$. 
Fix then a neighbourhood $U$ of $a$ on which $f$ is injective. 
Since $f^n$ converges to $a$ on $U$, 
$f^n(U)\subset U$ for all large enough integers $n\geq n_0$. 
Consider now a neighbourhood $U'$ of $a$ in $U$ whose images 
$U',f(U'),\dots,f^{n_0}(U')$ are all contained in $U$.
Such a set exists because $f$ is continuous and $a$ is a fixed point 
of $f$. By construction $f^n(U')\subset U$ for all $n\in \N$, and the restriction
of $f^n$ to $U'$ is injective as a composition of injective maps.

In the general setting, let us first check that in fact, the convergence of the subsequence 
$f^{n_k}$ to $a$ implies the convergence of the whole dynamics of an iterate $h=f^p$ to $a$. 
Pick again a small neighbourhood $U$ of $a$ in $\spc(b\om)$ and an integer $p=n_{k_0}$ such that 
$f^p(U)\subset U$. The map $h:=f^p$ restricts to $U$ to a local diffeomorphism from 
$U$ to itself, whose sequence of images $h^n(U)$ is obviously decreasing ({\it i.e.}
$h^i(U)\supset h^{i+1}(U)$). Observe then that the subsequence $(h^{n_k'})$ defined by
$n_k':=E(n_k/p)+1$ converges uniformly to $a$ on $U$. Indeed,
$h^{n_k'}=f^{pn_k'}=f^{n_k+i}$ with  $i<p$, so 
  $h^{n_k'}(U)\subset \cup_{i\leq p}f^i(f^{n_k}(U))$. Since $f^{n_k}(U)$ is close to 
$a$ by hypothesis (for $k$ large enough) and $a$ is a fixed point of $f$, the continuity 
of $f$ implies that $h^{n_k'}(U)$ is also close to $a$. Since the sequence $h^n(U)$ 
decreases, it thus converges to $a$. Replacing $f$ by $h$, we can therefore apply the above argument, 
so a \nbd of $a$ is indeed spherical.\cqfd

Consider now the situation when $f^{n_k}$ does not converge to $a$ on a \nbd of $a$.
 Let us first describe the figure and notation. As in the previous section, we assume 
 that  $\om$ is strongly convex in a \nbd $O$ of $a$, that $a$ is the origin 
and that $\om\cap O$ is contained in $\{\re z_1\geq 0\}$. We put  
$\om_\eps:=\om\cap O\cap\{\re z_1\leq\eps\}$,
$U_\eps:=b\om\cap O\cap\{\re z_1\leq\eps\}$ and we assume without loss of generality that 
$\om_1\Subset O$. Also since all the arguments to come are purely local and occur in $O$, 
we will consider in the sequel that $f$ extends smoothly to the boundary (lemma \ref{sext}), without explicitly mentionning any further the necessary restriction of $f$ to $O$.
 The non-convergence of $f^{n_k}$ means the existence 
of a sequence of points $z_i\in b\om$ tending to $a$, and integers $k_i$ such that 
the points $f^{n_{k_i}}(z_i)$ lay out of a fixed \nbd of $a$, say $U_1$. Since 
$a$ is fixed by $f^{n_{k_i}}$, we can even assume that 
$f^{n_{k_i}}(z_i)\in bU_1=b\om\cap O\cap\{\re z_1=1\}$ by moving $z_i$ closer to $a$. 
Finally, put $f_i:=f^{n_{k_i}}$ and define $\eps_i$ by $z_i\in\{\re z_1=\eps_i\}$.
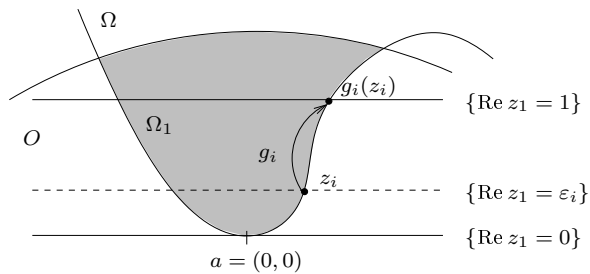
\begin{figure}[h]
\begin{center}
\input situation.pstex_t
\end{center}
\end{figure}

\noindent The main point of this section is that $f^{n_k}$ has a strong 
expanding behaviour.
\begin{prop}\label{dilat} (see also \cite{moi2})
For all $\eps$ there exists an integer $k=k(\eps)$ such that 
$f_k(U_\eps)\supset U_1\priv U_\eps$.
\end{prop}
\noindent The sphericity near $a$ is a direct consequence of this proposition :
\begin{cor}\label{hardspher}
If $(f^{n_k})$ does not converge to $a$ in a \nbd of $a$ then $b\om$
is spherical near $a$.
\end{cor}
\noindent{\it Proof :} Fix an open contractible set $V$ compact in 
$U_1\priv\{a\}$. For 
$\eps$ small enough, $V\subset U_1\priv U_\eps$ and there is an integer $k_\eps$
such that $f_{k_\eps}(U_\eps)\supset V$. Moreover, there are no critical 
value of ${f_{k_\eps}}_{|U_\eps}$ inside $V$ because both $U_\eps$ and $V$
are strictly pseudoconvex (see \cite{fornaess}). Since $V$ is simply connected,
there exists an inverse branch of ${f_{k_\eps}}_{|U_\eps}$ on $V$, which 
means a CR-diffeomorphism $h_\eps:V\lra U_\eps$ with $f_{k_\eps}\circ h_\eps=$Id.
The sequence $h_\eps$ is therefore contracting on $V$, and 
theorem \ref{wrlocal} implies that $V$ is spherical. We have thus proved the 
local sphericity of $U_1\priv\{a\}$, which even proves the sphericity of 
$U_1$ because $a$ is a strictly pseudoconvex point. Indeed, 
Chern-Moser's theory expresses the sphericity of an open  strictly pseudoconvex 
hypersurface by the vanishing of a continuous invariant tensor. 
Since this tensor vanishes on $U_1\priv\{a\}$, it also vanishes on the 
whole of $U_1$ so $U_1$ itself is spherical. In the spirit of \cite{moi2},
It would be pleasant to get a more down-to-earth argument for this last point.
\cqfd

The proof of proposition \ref{dilat} relies on the following lemma.
\begin{lemma}\label{hopf} For all $\eps$ there exists a diverging sequence $c_i\lra +\infty$
such that for all $p\in U_1$ with $f_i(p)\notin U_\eps$ we have :
$$
\Vert f_i'(p)u\Vert\geq c_i\Vert u\Vert \hspace{.5cm} \forall u\in T_p^\C b\om.
$$
\end{lemma} 
\noindent{\it Proof :} The idea is that Hopf's lemma gives estimates 
on the normal derivative of $f_i$, which transfer automatically 
to complex tangential estimates in strictly pseudoconvex geometry. 
For $p\in U_1$, let $\vec N(p)$ be the unit vector normal to $b\om$ 
pointing inside $\om$ and 
$$
B_\delta^+(p):=B(p+\delta\vec N(p),\delta)\cap
\{\langle \vec N(p),\cdot\rangle\geq \delta\}.
$$
When $\delta$ is small enough but fixed, $B_\delta^+(p)$ is in $\om$ and its image 
by $f_i$ for $i$ large is in $\om_{\frac{\eps}{2}}$ ($f_i$ converges to $a$ inside $\om$). Thus if $f_i(p)\notin U_\eps$,
the non-positive \psh\ function $\phi:=-\langle\vec N(f_i(p)),f_i(\cdot)-f_i(p)\rangle$ 
vanishes at $p$ while it is less
than $-c\eps^2$ on $B_\eps^+(p)$ ($c$ is a constant depending only on the 
curvature of $b\om$ at $a$). Hopf's lemma then asserts that 
$$
n_i(p):=\langle f_i'(p)\vec N(p),\vec N(f_i(p))\rangle=\Vert \grad \phi(p)\Vert
\geq \frac{c'\eps^2}{\delta}.
$$
Since $\delta$ was arbitrary, we could take it much smaller than $\eps^2$,
so that the radial escape rate $n_i(p)$ is large. 
To transfer this radial estimate on the derivatives of $f_i$ to complex
tangential ones, consider the Levi form $\L$ of $b\om$ defined by
$$
\L(p,u):=\langle[u,iu],i\vec N(p)\rangle,\hspace{0.5cm} u\in T^\C_pb\om,
$$
where $u$ stands for the vector in $T^\C_pb\om$ as well as any extension of it 
to a vector field of $T^\C b\om$. The smoothness and strict pseudoconvexity 
of $U_1$ implies the existence of geometric constants $c_1$, $c_2$ such that 
$$
c_1\Vert u\Vert^2\leq \L(p,u)\leq c_2\Vert u\Vert^2 \hspace{.5cm} \forall p\in U_1,
\; \forall u\in T^\C_pb\om.
$$
Easy computations show that :
$$
c_2 \Vert f'_i(p)u\Vert^2\geq \L(f_i(p),f_i'(p)u)=n_i(p)\L(p,u)\geq 
c_1 n_i(p)\Vert u\Vert^2.
$$
Since $n_i(p)$ is large when $i$ is, this serie of inequalities implies lemma 
\ref{hopf}.\cqfd
The previous lemma asserts that $f_i$ dilates the complex tangential directions
of $b\om$ if $f_i(p)$ is not close to $a$. The last observation we need to make in 
order to prove proposition \ref{dilat} is that this ``complex tangential dilation''
property implies a genuine dilation.

A path $\gamma$ in $b\om$ will be called a complex path if $\dot \gamma(t)\in 
T^\C_{\gamma(t)} b\om$ for all $t$. Its euclidean length will be denoted by 
$\ell(\gamma)$. For $x,y\in U_1$, define the CR-distance $d^{\text{\tiny CR}}(x,y)$
between $x$ and $y$ as the infimum of the lengths of complex paths joining $x$ to 
$y$. The point is that the strict pseudoconvexity condition means that the complex 
tangential distribution is contact so complex paths can join any two points. 
Even more, the open set $U_1\priv U_\eps$ is $d^{\text{\tiny CR}}$-bounded 
(see theorem 4 of \cite{nsw}, or \cite{moi1}).\\
\noindent{\it Proof of proposition \ref{dilat} :} Fix $\tau>0$ such that 
$B^{\text{\tiny CR}}(z_i,\tau)\subset U_\eps$ for all $i$ large enough. 
Since $U_1\priv U_\eps$ is $d^{\text{\tiny CR}}$-bounded, it is enough to prove that 
$$
bf_i(B^{\text{\tiny CR}}(z_i,\tau))\cap B^{\text{\tiny CR}}(f_i(z_i),c_i\tau)\cap 
\big(U_1\priv U_\eps\big) =\emptyset
$$
because $c_i\tau$ can be made greater than the CR-diameter of $U_1\priv U_\eps$.
Take a point $x\in bf_i(B^{\text{\tiny CR}}(z_i,\tau))\cap U_1\priv U_\eps$
and let us prove that 
\begin{equation}\label{eq1}
d^{\text{\tiny CR}}(f_i(z_i),x)\geq c_i\tau.
\end{equation}
Consider an arc-length parameterized complex path $\gamma$ in $U_1\priv U_\eps$
joining $f_i(z_i)$ to $x$. Since $f_i$ is a local CR-diffeomorphism  at each point 
of $B^{\text{\tiny CR}}(z_i,\tau)$ whose image lies in the strictly pseudoconvex part 
of $b\om$, the connected component of $f_i(z_i)$ in $\gamma\cap 
f_i(B^{\text{\tiny CR}}(z_i,\tau))$ can be lifted to a complex path 
$\wdt\gamma$ through $f_i$. Thus there exists $l\leq \ell(\gamma)$ and 
$\wdt \gamma:[0,l]\lra B^{\text{\tiny CR}}(z_i,\tau)$ joining $z_i$ to 
$bB^{\text{\tiny CR}}(z_i,\tau)$ such that $f_i\circ \wdt \gamma(t)=\gamma(t)$ for 
all $t\in [0,l]$. Since $\wdt\gamma(t)\in U_1$ and $f_i(\wdt\gamma(t))\in U_1\priv 
U_\eps$ for all $t$, the estimates obtained in lemma \ref{hopf} yield :
$$
\ell(\gamma) \geq  l  =  \int_0^l\Vert
 \dot{\gamma}(t)\Vert dt
 = \int_0^l \Vert f_i'(\wdt\gamma(t))\dot{\wdt\gamma}(t)\Vert dt \geq 
 c_i \int_0^l\Vert \dot{\wdt\gamma}(t)\Vert dt \\ 
  \geq 
c_i\ell(\wdt\gamma).
$$
This proves (\ref{eq1}) since $\wdt \gamma$ joins $z_i$ to 
$bB^{\text{\tiny CR}}(z_i,\tau)$ 
(so $\ell(\wdt\gamma)\geq \tau$) and  $\gamma$ is any complex path joining 
$f_i(z_i)$ to $x$.\cqfd

{ 
\footnotesize
\bibliography{bib.bib}
\bibliographystyle{abbrv}
}
\end{document}